\documentclass[12pt]{amsart}
\usepackage{amscd}
\def\NZQ{\Bbb}               
\def\NN{{\NZQ N}}
\def\QQ{{\NZQ Q}}
\def\ZZ{{\NZQ Z}}

\def\FFF{{\NZQ F}}

\def\frk{\frak}               

\def\pp{{\frk p}}

\def\qq{{\frk q}}

\def\mm{{\frk m}}

\def\opn#1#2{\def#1{\operatorname{#2}}} 
\opn\chara{char} \opn\length{\ell} \opn\pd{pd} \opn\rk{rk}
\opn\projdim{proj\,dim} \opn\injdim{inj\,dim} \opn\rank{rank}
\opn\depth{depth} \opn\codepth{codepth} \opn\grade{grade} \opn\height{height}
\opn\embdim{emb\,dim} \opn\codim{codim}

\opn\Tr{Tr} \opn\bigrank{big\,rank}
\opn\superheight{superheight}\opn\lcm{lcm}
\opn\trdeg{tr\,deg}%
\opn\reg{reg} \opn\lreg{lreg} \opn\skel{skel}
\opn\div{div} \opn\Div{Div} \opn\cl{cl} \opn\Cl{Cl}
\opn\Spec{Spec} \opn\Supp{Supp} \opn\supp{supp} \opn\Sing{Sing}
\opn\Ass{Ass}   \opn\grade{grade} \opn\cd{cd}
\opn\Ann{Ann} \opn\Rad{Rad} \opn\Soc{Soc}    \opn\Fitt{Fitt} \opn\Max{Max}
\opn\Sym{Sym} \opn\Ker{Ker} \opn\Coker{Coker} \opn\Im{Im}
\opn\Hom{Hom} \opn\Tor{Tor} \opn\Ext{Ext} \opn\End{End}
\opn\Aut{Aut} \opn\id{id} \opn\ini{in}

\opn\nat{nat}\opn\it{it}
\opn\pff{proof}
\opn\Pf{proof} \opn\GL{GL} \opn\SL{SL} \opn\mod{mod} \opn\ord{ord}
\opn\Hilb{Hilb}\opn\ara{ara} \opn\rdim{rdim} \opn \subsetneqq{subsetneqq}
\opn\aff{aff} \opn\con{conv} \opn\relint{relint} \opn\st{st}
\opn\lk{lk} \opn\cn{cn} \opn\core{core} \opn\vol{vol}
\opn\link{link} \opn\star{star} \opn\skel{skel}\opn\Kd{Krull-dim}
\opn\gr{gr}

\def\pot#1#2{#1[\kern-0.28ex[#2]\kern-0.28ex]}

\opn\dirlim{\underrightarrow{\lim}}
\opn\inivlim{\underleftarrow{\lim}}
\let\sect=\cap
\let\dirsum=\oplus
\let\tensor=\otimes
\let\iso=\cong
\let\Union=\bigcup

\let\Dirsum=\bigoplus

\let \varsubsetneq=\varsubsetneq
\def\Implies{\ifmmode\Longrightarrow \else
     \unskip${}\Longrightarrow{}$\ignorespaces\fi}
\def\implies{\ifmmode\Rightarrow \else
     \unskip${}\Rightarrow{}$\ignorespaces\fi}
\def\iff{\ifmmode\Longleftrightarrow \else
     \unskip${}\Longleftrightarrow{}$\ignorespaces\fi}

\let\:=\colon
\newtheorem{Theorem}{Theorem}[section]
\newtheorem{Lemma}[Theorem]{Lemma}
\newtheorem{Corollary}[Theorem]{Corollary}
\newtheorem{Proposition}[Theorem]{Proposition}
\newtheorem{Remark}[Theorem]{Remark}

\newtheorem{Example}[Theorem]{Example}

\newtheorem{Definition}[Theorem]{Definition}

\let\epsilon\varepsilon
\let\phi=\varphi
\let\kappa=\varkappa
\def\qed{\ifhmode\textqed\fi
   \ifmmode\ifinner\quad\qedsymbol\else\dispqed\fi\fi}
\def\textqed{\unskip\nobreak\penalty50
    \hskip2em\hbox{}\nobreak\hfil\qedsymbol
    \parfillskip=0pt \finalhyphendemerits=0}
\def\dispqed{\rlap{\qquad\qedsymbol}}

\opn\ini{in} \opn\inim{inm} \opn\rate{rate}

\opn\codim{codim}

\begin{document}

\title{Relative Cohen--Macaulayness of bigraded modules}

\author{Ahad Rahimi}

\subjclass[2000]{13C14, 13D45, 13D02, 16W50}

\address{Ahad Rahimi, Department of Mathematics, Razi University, Kermanshah, Iran} \email{ahad.rahimi@razi.ac.ir}
\maketitle
\address{}

 \maketitle
\begin{abstract}
In this paper we study the local cohomology of all  finitely generated
bigraded modules over a standard bigraded polynomial ring which have only one nonvanishing local cohomology with respect to one of the irrelevant bigraded ideals.
\end{abstract}
\maketitle
\section*{Introduction}
Let $S=K[x_1, \dots, x_m, y_1, \dots, y_n]$ be the standard bigraded polynomial ring over a field $K$ with bigraded
irrelevant ideals $P$ generated by all elements of degree $(1,0)$,
and $Q$ generated by all elements of degree $(0,1)$. In other words, $P=(x_1, \dots, x_m)$ and $Q=(y_1, \dots, y_n)$. Let $q\in \ZZ$ and $M$ be a finitely generated bigraded $S$-module such that $H^i_{Q}(M)=0$ for $i\neq q$. Thus $\grade(Q, M)=\cd(Q, M)=q$, where $\cd(Q, M)$ is the cohomological dimension of $M$ with respect to $Q$. Our aim is to characterize all finitely generated $S$-modules which have this property.
In this case,  we call $M$ to be relative Cohen--Macaulay with respect to $Q$ and call the common number $q$, relative dimension of $M$ with respect to $Q$ which will be denoted by $\rdim(Q, M)=q$.
We observe that ordinary Cohen--Macaulay modules are  special cases of our definition. In fact, if we assume that $P=0$, then $m=0$, and $Q=\mm$ is the unique graded maximal ideal of $S$. Therefore  $\depth M=\grade(\mm, M)=\cd(\mm, M)=\dim M$. We set $K[y]=K[y_1, \dots, y_n]$.
In Section 1 we show that $M$ is relative Cohen--Macaulay with respect to $Q$ of relative dimension $q$ if and only if $M_k=\Dirsum_jM_{(k,j)}$  is  a finitely generated Cohen--Macaulay $K[y]$-module of dimension $q$  for all $k$ and this is also equivalent to say that $M$ is Cohen--Macaulay as  $K[y]$-module of dimension $q$.
   In Section 2, we set  $H^q_{Q}(M)_j=\Dirsum_{k\in \ZZ}H^q_{Q}(M)_{(k,j)}$, and consider $H^q_{Q}(M)_j$ as a finitely generated graded $K[x]$-module where $K[x]=K[x_1, \dots, x_m]$. As a main result of this section,  we give a free resolution of $H^q_{Q}(M)_j$  which has length at most $m$. In particular, if $M$ is finitely generated bigraded  $S$-module of finite length, then each graded components $M_j$ has a free resolution of length $m$. By using this result we show that the regularity of $H^q_{Q}(M)_j$ is bounded for all $j$, i.e., there exists an integer $c$ such that  $-c\leq \reg H^q_{Q}(M)_j\leq c$, for all $j$.  The rest of this section is devoted to recall some algebraic properties of $H^q_{Q}(M)_j$.
   In Section 3,  we assume that $M$ be a Cohen--Macaulay $S$-module and show that  if $M$ is relative Cohen--Macaualy with respect to $Q$, then $M$ is relative Cohen--Macaulay with respect to $P$ and we have $\rdim (Q, M)+\rdim(P, M)=\dim M$. We deduce that, $M$ is relative Cohen--Macaualy with respect to $P$ or $Q$ if and only if $M$ satisfies  to this equation  $\dim (M/PM)+\dim(M/QM)=\dim M$. In the following section we assume that $M$ is relative Cohen--Macaulay with respect to $Q$. We prove the bigraded version of prime avoidance theorem which gives us a bihomogeneous  element $z\in Q$ for which $M/zM$ is relative Cohen--Macaulay with respect to $Q$ and relative dimension goes down by 1. As main result of this paper we prove that if $M$ is relative Cohen--Macaulay with respect to $Q$ then the following equality is always true:  $\rdim(Q,M)+\cd(P,M)=\dim M$. By using this result we show that if $M$ is relative Cohen--Macaulay with respect to $P$, as  well, then $M$ must be Cohen--Macaulay.  Some more applications of our main result are considered. Finally, we call $M$ to be maximal relative Cohen--Macaulay with respect to $Q$ if $M$ is relative Cohen--Macaulay with respect to $Q$ such that $\rdim (Q, M)=\dim K[y]=n$. We observe that maximal relative Cohen--Macaulay modules with respect to $Q$ are all finitely generated modules for which the sequence $y_1, \dots, y_n$ is an $M$-sequence, and free modules are the only modules which are maximal relative Cohen--Macaulay with respect to both irrelevant ideals.
\section{On the definition of relative Cohen--Macaulay }
Let $K$ be a field and  $S=K[x_1, \dots, x_m, y_1, \dots, y_n]$ be the standard bigraded polynomial ring. In other words, we set $\deg x_i=(1,0)$ and $\deg y_j=(0,1)$ for all $i,j$.
Let $M$ be a finitely generated bigraded $S$-module. We set $K[y]=K[y_1, \dots, y_n]$ and
\[
M_k=M_{(k,*)}=\Dirsum_{j\in \ZZ} M_{(k,j)}.
\]
 Then we view $M=\Dirsum_{k\in \ZZ}M_k$ as a graded module where each graded component $M_k$ itself is a finitely generated graded $K[y]$-module with grading $(M_k)_j=M_{(k,j)}$ for all $j$.
We shall use the following results \cite [Lemma 1.2.2]{RS2}.
 \begin{eqnarray}
\label{a}
 \dim_{K[y]}M=\sup\{\dim_{K[y]}M_k: k\in \ZZ\},
 \end{eqnarray}
 and
 \begin{eqnarray}
\label{b}
 \depth_{K[y]}M=\inf\{\depth_{K[y]}M_k: k\in \ZZ  \quad \text {with} \quad  M_k\neq 0\}.
 \end{eqnarray}
We may view  $M$ as graded $K[y]$-module which of course is not finitely generated in general. The module $M$ is called Cohen--Macaulay as $K[y]$-module if and only if $\depth_{K[y]} M=\dim_{K[y]}M$ where
\[
\depth_{K[y]} M=\sup\{n\in \NN:\quad \text{there is an $M$-sequence of length}\quad n\}.
\]
Let $R$ be a positively graded Noetherian ring and $M$ be a graded $R$-module.  We set $R_+=\dirsum_{i>0}R_i$ and denote by $\cd(R_+, M)$ the cohomological dimension of $M$ with respect to $R_+$ which is the largest integer $i$ for which  $H^i_{R_+}(M) \neq 0$.
When $M$ is a finitely generated $R$-module, by \cite[Lemma 3.4]{BrH} one has
\[
\cd(R_+,M)=\sup\{\dim_R M/\mm_0M: \mm_0\in \Max {R_0}\}.
\]
Thus if $(R_0,\mm_0)$ is a local ring, we have
\begin{eqnarray}
\label{formula 1}
 \cd(R_+,M)=\dim_R M/\mm_0M.
 \end{eqnarray}
We also denote by $\grade(R_+,M)$ the grade of $M$ with respect to $R_+$ which is the smallest integer $i$ for which  $H^i_{R_+}(M) \neq 0$. We also have the following characterization of $\grade(R_+, M)$,
\[
\grade(R_+,M)=\min\{i\in \NN_0:\Ext^i_R(R/R_+,M)\neq 0\}.
\]
Thus we see that $\grade(R_+,M)\leq \cd(R_+,M)\leq \dim M$ and if $(R, \mm)$ is a local ring we have $\grade(\mm,M)=\depth M$ and $\cd(\mm,M)=\dim M$.

In the following we give a necessary and sufficient condition for a finitely generated bigraded $S$-module $M$ which has only one nonvanishing local cohomology. This also generalizes \cite[Corollary 4.6]{RS}.
\begin{Proposition}
\label{one}
Let $M$ be a finitely generated bigraded $S$-module,  $q\in \ZZ$ and  $Q=(y_1, \dots, y_n)$. Then the following statements are equivalent:
\begin{itemize}
\item[{(a)}] $H^i_{Q}(M)=0$ for all $i\neq q$;
\item[{(b)}] $M_k=\Dirsum_jM_{(k,j)}$  is  finitely generated Cohen--Macaulay $K[y]$-module of dimension $q$  for all $k$;
\item[{(c)}] $M$ is Cohen--Macaulay as  $K[y]$-module of dimension $q$;
\item[{(d)}] $\grade(Q,M)=\cd(Q,M)$.
\end{itemize}
\end{Proposition}
\begin{proof}
(a)\implies (b), (c): We first observe that
\begin{eqnarray}
\label{a}
H^i_{Q}(M)=\Dirsum_k H^i_{Q}(M)_{(k,*)}=\Dirsum_k H^i_{Q}(M_{(k,*)}).
\end{eqnarray}
The second equality follows from the fact that $ H^i_{Q}(M)_{(k,*)}= H^i_{Q}(M_{(k,*)})$,
as can be seen from the definition of local cohomology using the \v{C}ech complex.
Now let $H^i_{Q}(M)=0$ for all $i\neq q$,  it follows that $H^i_{Q}(M_{(k,*)})=0$ for all $k$ and $i\neq q$. Hence $M_k=M_{(k,*)}$ is a finitely generated Cohen--Macaulay $K[y]$-module of dimension $q$  for all $k$. By using (1) and (2) we see that $M$ is Cohen--Macaulay as  $K[y]$-module of dimension $q$.

(b)\implies(a): In view of (4), it is clear.

(c)\implies (b): Let $M$ be a Cohen--Macaulay as  $K[y]$-module of dimension $q$. Then by using (1) and (2), we have
\[
\dim_{K[y]}M_k\geq \depth_{K[y]}M_k\geq \depth_{K[y]}M=q=\dim_{K[y]}M\geq \dim_{K[y]}M_k.
\]
Thus $\depth_{K[y]}M_k=\dim_{K[y]}M_k=q$, as required.

(a)\implies (d): is obvious.
\end{proof}

Now we can make the following definition:
\begin{Definition}{\em
Let $M$ be a finitely generated bigraded $S$-module and $q\in \ZZ$. We call $M$ to be relative Cohen--Macaulay with respect to $Q$ if and only if $M$ satisfies in one of the equivalent conditions of Proposition \ref{one}. For a relative Cohen--Macaulay module $M$ with respect to $Q$, we define the relative dimension $\rdim(Q, M)$ of $M$ to be the number $q$ in Proposition \ref{one}.}
\end{Definition}


\section{Algebraic invariants of the graded components of $H^q_{Q}(M)$ }

We set $K[x]=K[x_1, \dots, x_m]$. For all $j$, we set
\[
 H^q_{Q}(M)_j=\Dirsum_{k\in \ZZ}H^q_{Q}(M)_{(k,j)},
\]
  and consider $H^q_{Q}(M)_j$ as a finitely generated graded $K[x]$-module with grading $\big(H^q_{Q}(M)_j\big)_k=H^q_{Q}(M)_{(k,j)}$. In this section we are going to study the algebraic properties of  $H^q_{Q}(M)_j$. Let $F$ be a finitely generated bigraded free $S$-module. So that $ F= \Dirsum_{i=1}^{t}S(-a_{i},-b_{i})$. By using Formula 1 in \cite{AR1} we obtain
\begin{eqnarray}
\label{5}
H_{Q}^n(F)_j=\Dirsum_{i=1}^{t}\Dirsum_{\left|a\right|=-n-j+b_{i}}K[x](-a_{i})z^a.
\end{eqnarray}
 Thus, we may consider  $H_{Q}^n(F)_j$ as finitely generated graded free $K[x]$-module. With this observation we prove the following

\begin{Theorem}
\label{relative}
Let $M$ be a finitely generated bigraded  $S$-module which has a free resolution of the form
\[
\FFF: 0 \rightarrow F_{m+n}\overset {\phi_{m+n}} \longrightarrow F_{m+n-1} \rightarrow \cdots \rightarrow F_1 \overset {\phi_1}\longrightarrow F_0 \overset {\phi_0}\longrightarrow 0,
\]
where $ F_i= \Dirsum_{k=1}^{t_i}S(-a_{ik},-b_{ik})$. Let $M$ be relative Cohen--Macaulay with respect to $Q$ of relative dimension $q$. Then for all $j$, the $K[x]$-module $H^q_{Q}(M)_j$ has a free $K[x]$-resolution of length at most $m$ and of the form
\[
 0 \rightarrow H_{Q}^n(F_{m+n-q})_j \overset {\psi_{m+n-q}}\longrightarrow \cdots \rightarrow H_{Q}^n(F_{n-q+1})_j \overset {\psi_{n-q+1}}\longrightarrow \Ker \psi_{n-q} \rightarrow H^q_{Q}(M)_j \rightarrow 0,
\]
where the maps  $\psi_i:H_{Q}^n(F_{i})_j \rightarrow H_{Q}^n(F_{i-1})_j $ are induced by $\phi_i$ for all $i$.
\end{Theorem}
\begin{proof}
 As $M$ is relative Cohen--Macaulay with respect to $Q$. Proposition \ref {one} implies that $M$ is Cohen--Macaulay as $K[y]$-module of dimension $q$.  Note that $\depth_SM\geq \depth_{K[y]}M=q$. Hence  $\projdim_SM\leq m+n-q$ and so $F_i=0$ for $i>m+n-q$.
Applying the functor $H_{Q}^n(-)_j$ to the resolution $\FFF$ yields a graded complex of free $K[x]$- modules
\begin{eqnarray*}
\label{5}
H_{Q}^n(\FFF)_j:  0 \rightarrow H_{Q}^n(F_{m+n-q})_j \rightarrow  H_{Q}^n(F_{m+n-q-1})_j\rightarrow \cdots \overset {\psi_{1}}\longrightarrow  H_{Q}^n(F_0 )_j\overset {\psi_{0}}\longrightarrow  0.
\end{eqnarray*}
By \cite[Theorem 1.1]{AR1} we have the following graded isomorphisms of $K[x]$-modules
\[
H_{Q}^{n-i}(M)_j\iso H_i(H_{Q}^n(\FFF)_j).
\]
Since  $H^i_{Q}(M)=0$ for all $i\neq q$. It follows that
\begin{eqnarray}
\label{6}
H_i(H_{Q}^n(\FFF)_j)=\left\{
\begin{array}{cc}
H^q_{Q}(M)_j & \text{for $i=n-q$},\\
0 & \text{for  $i\neq n-q$.}
\end{array}
\right.
\end{eqnarray}
This means that we have only one homology throughout of the complex $H_{Q}^n(\FFF)_j$.
 Thus the complex $H_{Q}^n(\FFF)_j$  breaks to the following exact sequences of $K[x]$-modules
\[
 0 \rightarrow H_{Q}^n(F_{m+n-q})_j \rightarrow \cdots \rightarrow H_{Q}^n(F_{n-q+1})_j \rightarrow \Im \psi_{n-q+1}\rightarrow 0,
\]
and
\[
0\rightarrow  \Im \psi_{n-q+1}\rightarrow \Ker \psi_{n-q}\rightarrow H^q_{Q}(M)_j \rightarrow 0,
\]
where $H^q_{Q}(M)_j=\Ker \psi_{n-q}/\Im \psi_{n-q+1}$.
Combining these two exact sequences and obtain the following free resolution for $H^q_{Q}(M)_j$ which has length at most $m$
\[
 0 \rightarrow H_{Q}^n(F_{m+n-q})_j \rightarrow \cdots \rightarrow H_{Q}^n(F_{n-q+1})_j \rightarrow \Ker \psi_{n-q} \rightarrow H^q_{Q}(M)_j\rightarrow 0.
\]
To complete our proof we only need to show that $\Ker \psi_{n-q}$ is free. Here we distinguish two cases: We first assume that $q=n$. Thus $\Ker \psi_{0}=H_{Q}^n(F_{0})_j $ and so $H_{Q}^n(M)_j$ has a free resolution of the form
\begin{eqnarray}
 0 \rightarrow H_{Q}^n(F_{m})_j \rightarrow \cdots \rightarrow H_{Q}^n(F_{1})_j \rightarrow H_{Q}^n(F_{0})_j \rightarrow H^n_{Q}(M)_j\rightarrow 0.
\end{eqnarray}
Now suppose  that $0\leq q< n$. Thus we have the short exact sequence
\[
0 \rightarrow \Ker \psi_1 \rightarrow H_{Q}^n(F_1)_j \rightarrow H_{Q}^n(F_0 )_j\rightarrow 0
\]
 in which $H_{Q}^n(F_0 )_j$ and  $H_{Q}^n(F_1)_j$ are free $K[x]$-modules. Hence the sequence is split exact and we have that $\Ker \psi_1$ is a free $K[x]$-module. We also have the following split exact sequence
\[
0 \rightarrow \Ker \psi_2 \rightarrow H_{Q}^n(F_2)_j \rightarrow \Im \psi_2 \rightarrow 0
\]
 in which $H_{Q}^n(F_2)_j$ and $\Im \psi_2=\Ker \psi_1$ are free $K[x]$-modules. This follows that $\Ker \psi_2=\Im \psi_3$ is free. We proceed the same argument as above and see from the last split exact sequence
 \[
0 \rightarrow \Ker \psi_{n-q} \rightarrow H_{Q}^n(F_{n-q})_j \rightarrow \Im \psi_{n-q} \rightarrow 0
\]
that $\Ker \psi_{n-q}$ is free.
\end{proof}
As a consequence we give a free resolution for each graded components of any finitely generated bigraded $S$-module $M$ for which $\cd(Q, M)=0$. This also includes all finitely generated modules of finite length.
\begin{Remark}{\em
\label{gama}
 Let $I$ be a ideal of Noetherian ring $R$ and $M$ be a finitely generated $R$-module, then $\cd(I,M)=0$ if and only if $H^0_{I}(M)=M$. This results from the following isomorphisms
\[
H_I^i(M)=H_I^i\big(M/H_I^0(M)\big) \quad \text {for all} \quad i>0,
\]
(See \cite[Corollary 2.1.7]{BS}) and the fact that $M/H_I^0(M)$ is $I$-torsion free, i.e.,

$H_I^0\big(M/H_I^0(M)\big)=0$.  See \cite[Lemma 2.1.2]{BS}.}
\end{Remark}
\begin{Corollary}
Let $M$ be a finitely generated  bigraded $S$-module with $\cd(Q, M)=0$. Then for all $j$, the finitely generated graded $K[x]$-module $M_j$ has the resolution of the form
\[
 0 \rightarrow H_{Q}^n(F_{m+n})_j \rightarrow \cdots \rightarrow H_{Q}^n(F_{n+1})_j \rightarrow \Ker \psi_{n} \rightarrow M_j\rightarrow 0.
\]
\end{Corollary}
\begin{proof}
The assertion follows from Theorem \ref{relative} and Remark \ref{gama}.
\end{proof}


Let  $M$ be a finitely generated graded $K[x]$-module with graded minimal free resolution
\[
\FFF: 0 \rightarrow F_{k} \longrightarrow F_{k-1} \rightarrow \cdots \rightarrow F_1 \longrightarrow F_0 \rightarrow 0.
\]
The Castelnuovo--Mumford regularity of $M$ is the invariant
\[
\reg(M)=\max\{b_i(\FFF)-i: i\geq 0\}
\]
where $b_i(\FFF)$ denotes the maximal degree of the generators of $F_i$.  In the following we show that if $M$ be a relative Cohen--Macaulay  $S$-module with respect to $Q$ of relative dimension $q$, then the regularity of $H^q_{Q}(M)_j$ is bounded for all $j$.
\begin{Proposition}
Let $M$ be a finitely generated $S$-module which is relative Cohen--Macaulay with respect to $Q$ of relative dimension $q$. Then the function $f_M(j)=\reg H^q_{Q}(M)_j$ is bounded.
\end{Proposition}
\begin{proof}
 By \cite[Proposition 2.6]{AR1} and in view of (6),  the function $f_M$ is bounded below. Thus it suffices to show that $f_M$ is bounded above.  We first assume that $q=n$. By (7) the $K[x]$-module $H_{Q}^n(M)_j$ has a free resolution of the form
\[
 0 \rightarrow H_{Q}^n(F_{m})_j \rightarrow \cdots \rightarrow H_{Q}^n(F_{1})_j \rightarrow H_{Q}^n(F_{0})_j \rightarrow H^n_{Q}(M)_j\rightarrow 0,
\]
where $ F_i= \Dirsum_{k=1}^{t_i}S(-a_{ik},-b_{ik})$. By (5) we have
\[
H_{Q}^n(F_i)_j=\Dirsum_{k=1}^{t_i}\Dirsum_{\left|a\right|=-n-j+b_{ik}}K[x](-a_{ik})z^a.
\]
 Hence
\[
\reg H^q_{Q}(M)_j\leq \max \{b_i(H_{Q}^n(\FFF)_j)-i:i\geq 0\},
\]
where $b_i(H_{Q}^n(\FFF)_j)$ is the maximal degree of the generators of $H_{Q}^n(F_i)_j$. Thus we conclude that
\[
\reg H^q_{Q}(M)_j\leq \max_{i, k}\{a_{ik}-i\}=c.
\]
Now let $0\leq q<n$, and consider the following exact sequences of $K[x]$-modules which we observed in the proof of Theorem \ref{relative}.
\begin{eqnarray}
\label{a}
 0 \rightarrow H_{Q}^n(F_{m+n-q})_j \rightarrow \cdots \rightarrow H_{Q}^n(F_{n-q+1})_j \rightarrow \Im \psi_{n-q+1}\rightarrow 0,
\end{eqnarray}

and
\begin{eqnarray}
\label{a}
0\rightarrow  \Im \psi_{n-q+1}\rightarrow \Ker \psi_{n-q}\rightarrow H^q_{Q}(M)_j \rightarrow 0,
\end{eqnarray}
where $H^q_{Q}(M)_j=\Ker \psi_{n-q}/\Im \psi_{n-q+1}$. Thus (8) yields
\begin{eqnarray*}
\reg \Im \psi_{n-q+1} &\leq & \max \{b_i(H_{Q}^n(\FFF)_j)+n-q+1-i:i\geq n-q+1\}\\
                   & =& \max_{i\geq n-q+1,k}\{a_{ik}+n-q+1-i\}=c,
\end{eqnarray*}
for some number $c$, and by (9) we have
\[
\reg H^q_{Q}(M)_j\leq \max \{\reg \Im \psi_{n-q+1}-1,  \reg \Ker \psi_{n-q}\}.
\]
Since $\reg \Im \psi_{n-q+1}$ is bounded above, to complete our prove, it sufficient to show that   $\reg \Ker \psi_{n-q}$ is bounded above.
To do so, from the short exact sequence
\[
0 \rightarrow \Ker \psi_1 \rightarrow H_{Q}^n(F_1)_j \rightarrow H_{Q}^n(F_0 )_j\rightarrow 0,
\]
it follows that
\[
\reg \Ker \psi_1 \leq \max\{\reg H_{Q}^n(F_1)_j,  \reg H_{Q}^n(F_0 )_j+1\}.
\]
Observing that $\reg H_{Q}^n(F_0 )_j=\max_k\{a_{0k}\}=a$ and $\reg H_{Q}^n(F_1 )_j=\max_k\{a_{1k}\}=b$ for some numbers $a$ and $b$. Thus the regularity of $\Ker \psi_1$ is bounded above.
As $\Ker \psi_1=\Im \psi_2$, the exact sequence
\[
0 \rightarrow \Ker \psi_2 \rightarrow H_{Q}^n(F_2)_j \rightarrow \Im \psi_2 \rightarrow 0
\]
 yields that $\reg \Ker \psi_2$ is bounded above. We proceed the same argument as above and observe from the last exact sequence
 \[
0 \rightarrow \Ker \psi_{n-q} \rightarrow H_{Q}^n(F_{n-q})_j \rightarrow \Im \psi_{n-q} \rightarrow 0
\]
that $\reg \Ker \psi_{n-q}$ is bounded above and therefore the regularity of $\reg H^q_{Q}(M)_j$ is bounded above for all $j$, as required.
\end{proof}
In the remaining of this section we are going to recall some known results about $H_{Q}^q(M)$.
Let $R$ be a graded ring and $N$ a graded $R$-module. The $R$-module $N$ is called tame, if there exists an integer
$j_0$ such that either
\[
N_j=0\quad\text{for all}\quad j\leq j_0, \quad \text{or}\quad
N_j\neq 0\quad \text{for all}\quad j\leq j_0.
\]
\begin{Proposition}
\label{properties} Let $M$ be a finitely generated bigraded $S$-module which is relative Cohen--Macaulay  with respect to $Q$ of  relative dimension $q$. Then the following statements hold:
\begin{itemize}

\item[{(a)}] $\Ass_{K[x]} H_{Q}^q(M)_j$ is asymptotically stable for $j\longrightarrow -\infty$, i.e., there exists an integer $j_0$ such that $\Ass_{K[x]} H_{Q}^q(M)_j=\Ass_{K[x]} H_{Q}^q(M)_{j_0}$ for all $j\leq j_0$;
    \item[{(b)}] $\Ass_
    S H_{Q}^q(M)<\infty$;
 \item[{(c)}]   $H_{Q}^q(M)$ is tame.

\end{itemize}
\end{Proposition}
\begin{proof}
Part (a) follows from \cite[Lemma 5.4]{BrH} simply because $H_{Q}^i(M)$ is finitely generated for $i<q$. Part (b) follows from  \cite[Remark 5.5(B)]{BrH} and part (c) results from (a).
\end{proof}

Let $R$ be a Noetherian local ring and $I$ be an ideal of $R$. An $R$-module $N$ is said to be $I$-cofinite if
$\Supp N \subseteq V(I)$ and $\Ext_R^i(R/I,N)$ is finitely generated for all $i\geq 0$.
\begin{Proposition}
Let $M$ be a finitely generated bigraded $S$-module which is relative Cohen--Macaulay  with respect to $Q$  of relative dimension  $q >0$. Then $H_{Q}^q(M)$ is not finitely generated, but it is $Q$-cofinite.
\end{Proposition}
\begin{proof}
By Proposition \ref{properties} part (c), $H_{Q}^q(M)$ is tame. In fact,  $H_{Q}^q(M)_j\neq 0$ for $j\ll 0$. (See \cite[Remark 2.6]{RS}).
So $H_{Q}^q(M)$ can not be finitely generated, because otherwise $H_{Q}^q(M)_j = 0$ for $j\ll 0$. Now we show that $H_{Q}^q(M)$ is $Q$-cofinite. As $\Supp H_{Q}^q(M) \subseteq V(Q)$, thus we only need to show that $\Ext_S^i(S/Q,H_{Q}^q(M))$ is finitely generated for all $i\geq 0$. We consider the following spectral sequence of $S$-modules
\[
\Ext_S^i(S/Q, H_{Q}^j(M))\underset{i}
\Longrightarrow \Ext_S^{i+j}(S/Q, M).
\]
Since $H_{Q}^j(M)=0$ for all $j\neq q$. Thus we have the following isomorphisms of $S$-modules
$\Ext_S^i(S/Q, H_{Q}^q(M))\iso \Ext_S^{i+q}(S/Q,M)$ which are always finitely generated for all $i\geq 0$.
\end{proof}
\begin{Remark}{\em
Let $R$ be a graded Noetherian  ring and $M$ be a finitely generated graded $R$-module. The finiteness dimension of $M$ with respect to $R_+$ is defined as
\[
f(M)=\inf\{i\in \NN_0: H_{R_+}^i(M) \quad \text {is not finitely generated} \}.
\]
The $R_0$-module $H_{R_+}^i(M)_j$ are finitely generated and $H_{R_+}^i(M)_j=0$ for $j\gg 0$. Thus we may write
\[
f(M)=\inf\{i\in \NN_0: H_{R_+}^i(M)_j\neq 0  \quad \text {for infinitely many $j<0$ } \}.
\]
The finite length dimension of $M$ is defined as
\[
g(M)=\inf\{i\in \NN_0: \length \big( H_{R_+}^i(M)_j\big)=\infty \quad \text {for infinitely many $j<0$ } \}.
\]
Thus we observe that $f(M)\leq g(M)$.
Now we assume that $M$ is relative Cohen--Macaulay with respect to $Q$ of relative dimension $q>0$. Hence  $f(M)=g(M)=q$. Then by \cite[Theorem 4.11 (a)]{BRS}, there exists a polynomial $\widetilde{Q}\in \QQ[x]$ such that the multiplicity
\[
e(H_{Q}^q(M)_j )=\widetilde{Q} (j) \quad \text {for } \quad j\ll0.
\]
By \cite[Corollary 2.4]{B} the minimal number of generators of $H_{Q}^q(M)_j$ is also a polynomial for $j\ll 0$. More precisely, there exists a polynomial  $\widetilde{Q}\in \QQ[x]$ of degree  at most $q$ such that
\[
\mu(H_{Q}^q(M)_j )=\widetilde{Q} (j) \quad \text {for } \quad j\ll0.
\]}
\end{Remark}
\section{Cohen--Macaulay modules which are relative Cohen--Macaulay}

Let $(R,\mm)$ be a Noetherian local ring, and $M$ a $R$-module. Then
\begin{eqnarray}
\label{formula 4}
\grade (I,M)\leq \dim M-\dim M/IM \quad \text {for all ideals} \quad I\subseteq \mm,
\end{eqnarray}
and equality holds if $M$ is Cohen--Macaulay (see \cite [Theorem 2.1.2]{BH}). So for a Cohen--Macaulay module the grade an
arbitrary ideal is given by its codimension. Now we are going to give a very explicit criteria for all finitely generated bigraded Cohen--Macaulay modules which are relative Cohen--Macaulay with respect to one of the irrelevant bigraded ideals $P$ and $Q$.
In fact, not all Cohen--Macaulay $S$-modules are relative Cohen--Macaulay. Obvious examples are hypersurface rings which are Cohen--Macaulay but have two nonvanishing local cohomology.

In the sequel, the graded version of (10) will be considered.

\begin{Proposition}
\label{rcm}
Let $M$ be a finitely generated bigraded Cohen--Macaulay $S$-module. Then the following statements are equivalent:
\begin{enumerate}
\item[{(a)}] $M$ is relative Cohen--Macaulay with respect to $P$;
\item[{(b)}] $M$ is relative Cohen--Macaulay  with respect to $Q$;
\item[{(c)}] $M$ is relative Cohen--Macaulay  with respect to $P$ and $Q$ with $$\rdim (P, M)+\rdim(Q, M)=\dim M.$$
\end{enumerate}
\end{Proposition}
\begin{proof}
(a)\implies (b), (c):
By (10) we have $\grade (Q, M)=\dim M-\dim M/QM$, since $M$ is Cohen--Macaulay.
 In view of (3)  we have $\grade (P,M)= \cd(P, M)=\dim M/QM$ and $\cd(Q, M)=\dim M/PM$. By using these facts and again (10)  we have
\begin{eqnarray*}
\grade (Q, M)&=&\dim M-\grade(P, M)\\
                                 &=& \dim M-(\dim M -\dim M/PM)\\
                                 &=& \dim M/PM=\cd(Q, M).
\end{eqnarray*}
Thus we conclude that $M$ is relative Cohen--Macaulay  with respect to $Q$ and that $\rdim (P, M)+\rdim(Q, M)=\dim M$.

(b) \implies (a), (c):  is proved the same way, and (c) \implies (a), (b) is clear.
 \end{proof}
 \begin{Corollary}
Let $M$ be a finitely generated bigraded Cohen--Macaulay $S$-module. Then,  $M$ is relative Cohen--Macaulay with respect to $P$ or $Q$ if and only if
$$\dim (M/QM)+\dim(M/PM)=\dim M.$$
 \end{Corollary}
 \begin{proof}
 The assertion follows from Proposition \ref{rcm} and (10).
 \end{proof}
 In \cite {HR}  for a given bigraded $S$-module $M$ we define the bigraded
Matlis-dual of $M$ to be  $M^\vee$  where the $(i,j)$th bigraded
component of $M^\vee$ is given by $\Hom_K(M_{(-i,-j)},K)$. There exists a convergent spectral
sequence
\[
E^2_{i,j}=H^{m-j}_{P}\big(H^i_{\mm}(M)^\vee \big)\underset{j}
\Longrightarrow H_{Q}^{i+j-m}(M)^\vee
\]
of bigraded $S$-modules where $m$ is the number of variables of $P$ and $\mm=P+Q$ is the unique graded maximal
ideal of $S$. The above spectral sequence degenerates when $M$ is
Cohen--Macaulay and one obtains  for all $k$ the following
isomorphisms of bigraded $S$-modules
\begin{eqnarray}
\label{formula }
 H^{k}_P \big(H^s_{\mm}(M)^\vee\big)\iso H^{s-k}_{Q}(M)^\vee
\end{eqnarray}
 where $s=\dim M$, see \cite[Corollary 2.6]{HR}.
 \begin{Corollary}
 Let $M$ be a finitely generated bigraded Cohen--Macaulay  $S$-module which is relative Cohen--Macaulay  with respect to $Q$. Then $H^s_{\mm}(M)^\vee$ is relative Cohen--Macaulay  with respect to $P$ and $Q$,  and we have
 \[
 \rdim(P, M)=\rdim (P, H^s_{\mm}(M)^\vee) \quad \text {and} \quad \rdim(Q, M)=\rdim (Q, H^s_{\mm}(M)^\vee).
 \]
 \end{Corollary}
  \begin{proof}

  By Proposition \ref{rcm} and (11) we have
 \begin{eqnarray*}
 \rdim(P, M)=\dim M-\rdim(Q, M)& = & \dim M-(\dim M-\rdim (P, H^s_{\mm}(M)^\vee))\\
                                 & = & \rdim (P, H^s_{\mm}(M)^\vee).
 \end{eqnarray*}
 Since $M$ is Cohen--Macaulay, it follows that $H^s_{\mm}(M)^\vee$ is also Cohen--Macaulay with the same dimension as $M$. Using this fact and the first part together with Proposition \ref{rcm} we have
  \begin{eqnarray*}
 \rdim(Q, M)=\dim M-\rdim(P, M)& = & \dim H^s_{\mm}(M)^\vee-\rdim (P, H^s_{\mm}(M)^\vee))\\
                                 & = & \rdim (Q, H^s_{\mm}(M)^\vee),
 \end{eqnarray*}
 as desired.
 \end{proof}
 In the following  we are going to give some examples which show that Proposition \ref{rcm} can fail if $M$ is not Cohen--Macaulay. We first prove the following lemma.
 \begin{Lemma}
 \label{example}
  Let $I$ and $J$ be ideals  of $K[x]$ and $K[y]$, respectively. We set $R_0=K[x]/I$,  $R_1=K[y]/J$ and $R=R_0\tensor_KR_1$. Then the following statements hold:
 \begin{itemize}
\item[{(a)}] $R$ is relative Cohen--Macaulay with respect to $Q$ of relative dimension $q$ if and only if $R_1$ is Cohen--Macaulay of dimension $q$;
\item[{(b)}] $R$ is relative Cohen--Macaulay with respect to $P$ of relative dimension $p$ if and only if $R_0$ is Cohen--Macaulay of dimension $p$.
\end{itemize}
\begin{proof}
 In order to prove (a) we first observe that $R\iso K[x, y]/(I, J)$.  We consider the following isomorphisms of $S$-modules
\begin{eqnarray*}
H_{Q}^i(R)=H^i_{QR}(R)& \iso & \underset{k\geq 0}{\dirlim}\Ext_R^i\big(R/(QR)^k,R\big)\\
                         & \iso &   R_0\otimes_ K\underset{k\geq 0}{\dirlim}\Ext_{R_1}^i\big(R_1/(QR_1)^k,R_1\big)\\
           & \iso &   R_0\tensor_ K H_{QR_1}^i (R_1)\iso R_0\tensor_ K H_{Q}^i (R_1).
\end{eqnarray*}
Now the assertion follows from this observation.
Part (b) is proved the same way.
\end{proof}
 \end{Lemma}
 As a special case of Lemma \ref{example} we have the following example.
 \begin{Example}{\em
 \label{1}
 We consider the following standard bigraded ring
 \[
 R=\frac{K[x_1, x_2, y_1, y_2]}{(x_1^2, x_1x_2)}.
 \]
 One has $\depth R=2$ and $\dim R=3$,  and so $R$ is not Cohen--Macaulay. We set
 $R_0=K[x_1, x_2]/(x_1^2, x_1x_2)$ and $R_1=K[y_1, y_2]$.
  Since $R_1$ is Cohen--Macaualy of dimension $2$,  Lemma \ref{example} yields that $R$ is relative Cohen--Macaulay with respect to $Q$ of relative dimension $2$.  On the other hand,
 $\Ass(I)=\{(x_1),(x_1, x_2)\}$ where $I=(x_1^2, x_1x_2)$. Thus \cite[Proposition 9.1.4(a)]{BH} yields  $\grade(P, R)=0$. One has
 $\cd(P, R)=\dim R/QR=1$. So $R$ is not relative Cohen--Macaulay with respect to $P$.
   }
 \end{Example}
 \begin{Example}{\em
 \label{2}
 We consider the following standard bigraded ring
 \[
 R=\frac{K[x_1, \dots, x_m, y_1]}{(x_1y_1, \dots, x_my_1, y_1^2)}.
 \]
 We observe that $\depth R=0$ and $\dim R=m$,  and so $R$ is not Cohen--Macaulay. We see that $\grade(Q, R)\leq \grade(\mm, R)=\depth R=0$, and so $\grade(Q, R)=0$. One has $\cd(Q, R)=\dim R/PR=0$. Thus $R$ is relative Cohen--Macaulay with respect to $Q$. On the other hand,  $\grade(P, R)\leq \depth R=0$, and so $\grade(P, R)=0$. One has $\cd(P, R)=\dim R/QR=m$. Thus $R$ is not relative Cohen--Macaulay with respect to $P$.
  }
 \end{Example}

 \section{Relative Cohen--Macaulay with respect to $P$ and $Q$ must be Cohen--Macaulay }

 Let $M$ be a relative Cohen--Macaulay module with respect to $Q$. By using (3) and (10) we have
 \[
 \rdim (Q,M)+\cd(P,M)\leq \dim M.
 \]
 Our main result is to prove that equality holds in general. We shall need to use the following bigraded version of prime avoidance theorem.
 \begin{Lemma}
 \label{avoidance}
 Let ${\pp}_1,\dots, {\pp}_r$ be prime ideals of $S$ such that $Q\not \subseteq {\pp}_i$ for $i=1,\dots,r$ and $|K|=\infty$. Then there exists a bihomogeneous element $z\in Q$ of degree $(0,1)$ such that $z\not \in {\pp}_1\cup \dots \cup {\pp}_r$.
 \end{Lemma}
 \begin{proof}
Let $V$ be the $K$-vector
space spanned by $y_1, \dots, y_n$. Since $Q\not\subseteq p_i$ for $i=1, \dots, r$, it follows that
$V_i=V\sect p_i$ is a proper linear subspace of $V$. Since $|K|=\infty$, the
vector space $V$ can not be the finite union of proper linear subspaces.
Therefore, there exists $z\in V\setminus \Union_{i=1}^rV_i$. This is the
desired element $z$ of degree $(0,1$).
\end{proof}
\begin{Lemma}
\label{regular}
Let $M$ be a finitely generated bigraded $S$-module which is relative Cohen--Macaulay with respect to $Q$ with $\rdim(Q,M)>0$ and $|K|=\infty$. Then there exists a bihomogeneous  $M$-regular element $z\in Q$ such that $M/zM$ is relative Cohen--Macaulay with respect to $Q$ and we have
\[
\rdim(Q,M/zM)=\rdim(Q,M)-1.
\]
\end{Lemma}
\begin{proof}
By our assumptions we have $$\cd(Q,M)=\dim (M/PM)>0.$$ Let $\{{\pp}_1, \dots, {\pp}_r\}$ be the minimal prime ideals of $\Supp (M/PM)$. We claim that $Q\not \subseteq {\pp}_i$ for $i=1,\dots,r$. Assume that $Q \subseteq {\pp}_i$ for some $i$. Since  $P \subseteq {\pp}_i$ for all $i$. It follows that ${\pp}_i=P+Q=\mm$, and so $\dim (M/PM)=0$, a contradiction.
Since $\grade(Q,M)>0$, by \cite [Proposition 9.1.4 (a)]{BH} it follows that  $\Ass M \cap V(Q)=\emptyset$. Hence  $Q\not \subseteq \bigcup_{\qq\in \Ass M}\qq=Z(M)$ where $ Z(M)$ denotes the set zero divisors of $M$. By Lemma \ref {avoidance} we may choose a bihomogenous element $z\in Q$ which does not belong to any associated prime ideal of $M$ and not to any  minimal prime ideal of $\Supp (M/PM)$. In particular, this element is $M$-regular. By \cite [Proposition 9.1.2 (a)]{BH}  we have $\grade(Q,M/zM)=\grade(Q,M)-1$. To complete our prove, we only need to show that $\cd(Q,M/zM)= \cd(Q,M)-1$.  Since  $\grade(Q,M)=\cd(Q,M)$ and $\grade(Q,M/zM)\leq \cd(Q,M/zM)$ it follows that $\cd(Q,M)-1\leq \cd(Q,M/zM)$. By \cite[Corollary 2.3(i)]{DNT} the short exact sequence
 $
 0 \rightarrow zM \rightarrow M \rightarrow M/zM \rightarrow 0
 $
 yields $\cd(Q,M/zM)\leq \cd(Q,M)$. Thus we conclude that $\cd(Q,M)-1\leq \cd(Q,M/zM)\leq \cd(Q,M)$.
 Since the bihomogeneous element $z$ lies outside every minimal prime ideal of $\Supp (M/PM)$. It follows that
 \begin{eqnarray*}
 \cd(Q,M/zM) & = & \dim (M/zM)/P(M/zM)\\
               & = & \dim (M/PM)/z(M/PM)\\
               & < & \dim M/PM=\cd(Q,M),
 \end{eqnarray*}
 and so $\cd(Q,M/zM)=\cd(Q,M)-1$, as required.
\end{proof}
As main result of this section we prove the following
\begin{Theorem}
\label{rdim1}
Let $M$ be a finitely generated bigraded $S$-module which is relative Cohen--Macaulay with respect to $Q$ and $|K|=\infty$. Then we have
\[
\rdim(Q,M)+\cd(P,M)=\dim M.
\]
If in addition $M$ is relative Cohen--Macaulay with respect to $P$. Then we have
\[
\rdim(Q,M)+\rdim(P,M)=\dim M.
\]
\end{Theorem}
\begin{proof}
We prove the statement by induction on $\rdim(Q,M)$. Assume that $\rdim(Q,M)=0$. We consider the following spectral sequence
\[
H^i_{P}\big (H^j_{Q}(M)\big )\underset{i} \Longrightarrow H^{i+j}_{\mm}(M),
\]
where $\mm=P+Q $. As $H^j_{Q}(M)=0$ for  $j\neq 0$, then the above spectral sequence degenerates and one obtains  for all $i $ the following
isomorphism of bigraded $S$-modules,
\[
H^i_{P}\big (H^0_{Q}(M)\big) \iso H^{i}_{\mm}(M).
\]
Remark \ref{gama} implies that $H^0_{Q}(M)=M$,  and so  $H^i_{P}(M) \iso H^{i}_{\mm}(M)$.
Since $H^{\dim M}_{\mm}(M) \neq 0$, it follows that $H^{\dim M}_{P}(M)\neq 0$. Hence $\cd(P,M) \geq \dim M$.
 We also have that  $\cd(P,M) \leq \dim M$. Therefore $\cd(P,M)=\dim M.$
 Now suppose that $\rdim(Q,M)> 0$, and our desired equality has been proved for all finitely generated bigraded $S$-module $N$ such that $\rdim(Q,N)<\rdim(Q,M)$. We want to proof it for $M$. Since $\rdim(Q,M)> 0$, by Lemma \ref {regular} there exists a bihomogeneous  $M$-regular element $z\in Q$ such that $M/zM$ is relative Cohen--Macaulay with respect to $Q$ with  $\rdim(Q,M/zM)=\rdim(Q,M)-1$. Thus the induction hypothesis implies that
 \[
\rdim(Q,M/zM)+\cd(P,M/zM)=\dim M/zM.
\]
Since $z \in Q$, it follows that
\begin{eqnarray*}
\cd(P, M/zM) & = & \dim (M/zM)/Q(M/zM) \\
              & = & \dim M/(Q + (z))M\\
              & = & \dim M/QM=\cd(P, M).
\end{eqnarray*}
We also have $\dim M/zM=\dim M -1$. Therefore the desired equality follows.
\end{proof}
In Proposition \ref {rcm} we have shown that if $M$ is a Cohen--Macaulay module which is relative Cohen--Macaulay with respect to $Q$, then $M$ is relative Cohen--Macaulay with respect to $P$. Now we are going to show that if $M$ be any finitely generated bigraded $S$-module which is relative Cohen--Macaulay with respect to both $P$ and $Q$, then $M$ must be Cohen--Macaulay.
\begin{Corollary}
\label{Cohen--Macaulay}
Let $M$ be a finitely generated bigraded $S$-module which is relative Cohen--Macaulay with respect to both $P$ and $Q$ and $|K|=\infty$. Then $M$ is Cohen--Macaulay.
\end{Corollary}
\begin{proof}
We prove the corollary by induction on $\rdim(Q,M)$. We first assume that $\rdim(Q,M)=0$.  Theorem \ref{rdim1} implies that $\rdim(P,M)=\dim M$. Hence $\dim M = \grade(P,M)\leq \grade(\mm,M)=\depth M$ and so $M$ is Cohen--Macaulay.
Now we assume that $\rdim(Q,M)> 0$ and the desired result has been proved for all finitely generated bigraded $S$-module $N$ such that $\rdim(Q,N)<\rdim(Q,M)$. We want to proof it for $M$. As $\rdim(Q,M)> 0$, by Lemma \ref {regular} there exists a bihomogenous  $M$-regular element $z \in Q$ such that  $M/zM$ is relative Cohen--Macaulay with respect to $Q$ with  $\rdim(Q,M/zM)=\rdim(Q,M)-1$.   Applying our induction hypothesis implies that $M/zM$ is Cohen--Macaulay,  and so  $M$ is Cohen--Macaulay, as desired.
\end{proof}
\begin{Proposition}
Let $M$ be a finitely generated bigraded $S$-module which is relative Cohen--Macaulay with respect to $P$ and $Q$ of relative dimensions $p$ and $q$, respectively and $|K|=\infty$. Then the modules $H^i_{P}\big(H^q_{Q}(M)\big)$ and $H^i_{Q}\big(H^p_{P}(M)\big)$ are Artinian modules for all $i$ and we have the following isomorphism of bigraded module
\[
H^p_{P}\big(H^q_{Q}(M)\big)\iso H^q_{Q}\big(H^p_{P}(M)\big).
\]
\end{Proposition}
\begin{proof}
We consider the following spectral sequence
\[
H^i_{P}\big (H^j_{Q}(M)\big )\underset{i} \Longrightarrow H^{i+j}_{\mm}(M),
\]
where $\mm=P+Q$. As $M$ is relative Cohen--Macaulay with respect to $Q$ of  relative dimension $q$, we have the following isomorphisms of bigraded $S$-modules
\[
H^i_{P}\big (H^q_{Q}(M)\big )\iso H^{i+q}_{\mm}(M).
\]
   The modules $H^{i+q}_{\mm}(M)$ are Artinian for all $i$, and so the modules $H^i_{P}\big (H^q_{Q}(M)\big )$ are  Artinian for all $i$. Applying the isomorphism with $i=p$ yields
\begin{eqnarray}
\label{formula }
H^p_{P}\big (H^q_{Q}(M)\big )\iso H^{p+q}_{\mm}(M)=H^{\dim M}_{\mm}(M)\neq 0.
\end{eqnarray}
The last equality follows from Theorem \ref {rdim1}.
By the similar argument as above we see that the modules  $H^i_{Q}\big (H^p_{P}(M)\big )$ are Artinian for all $i$ and
 \begin{eqnarray}
\label{formula }
H^q_{Q}\big (H^p_{P}(M)\big )\iso H^{p+q}_{\mm}(M)=H^{\dim M}_{\mm}(M)\neq 0.
\end{eqnarray}
The desired isomorphism follows from (12) and (13).
\end{proof}
\begin{Lemma}
\label{Mayer}
Let $(R, \mm)$ be a Noetherian local ring, and $I$, $J$ ideals of $R$ such that $I+J=\mm$. Let $M$ be a finitely generated $R$-module with  $\cd(I, M)>0$,  $\cd(J, M)>0$ and $\cd(I,M)+\cd(J, M)=\dim M$. Then
\[
\cd(I\cap J, M)=\dim M-1.
\]
\end{Lemma}
\begin{proof}
The Mayer-Vietoris Sequence provides, the long exact sequence of $R$-modules
\[
\cdots \rightarrow H^{i-1}_{I\cap J}(M)\rightarrow H^i_{\mm}(M)\rightarrow H^i_{I}(M)\dirsum H^i_{J}(M)\rightarrow H^{i}_{I\cap J}(M)\rightarrow \cdots.
\]
Applying the long exact sequence with $i=\dim M$. Our assumptions yields that $0<\cd(I, M)<\dim M$ and $0<\cd(J, M)<\dim M$,  and so  $H^{\dim M}_{I}(M)\dirsum H^{\dim M}_{J}(M)=0$. Therefore we get the following exact sequence of $R$-modules \[
\cdots \rightarrow H^{\dim M-1}_{I}(M)\dirsum H^{\dim M-1}_{J}(M)\rightarrow H^{\dim M-1}_{I\cap J}(M)\rightarrow H^{\dim M}_{\mm}(M)\rightarrow 0 .
\]
Since $H^{\dim M}_{\mm}(M)\neq 0$, it follows that $H^{\dim M-1}_{I\cap J}(M)\neq 0$. Hence $\cd(I\cap J, M)\geq \dim M-1$. Thus we conclude that  $\dim M-1 \leq \cd(I\cap J, M)\leq \dim M$. The equality $\cd(I\cap J, M)= \dim M$ can not be the case, because by putting $i=\dim M$ in the above exact sequence yields $H^{\dim M}_{I\cap J}(M)=0$ .
Therefore $\cd(I\cap J, M)=\dim M-1$, as desired.
\end{proof}
\begin{Corollary}
\label{intersection}
Let $M$ be a finitely generated bigraded $S$-module which is relative Cohen--Macaulay with respect to $Q$ with relative dimension $\rdim(Q, M)>0$. Assume that $\cd(P, M)>0$ and $|K|=\infty$. Then we have
\[
\cd(P \cap Q, M)=\dim M-1.
\]
\end{Corollary}
\begin{proof}
The assertion follows from Theorem \ref{rdim1} and the graded version of Lemma \ref{Mayer}.
\end{proof}
\begin{Corollary}
Let the assumption be as in Corollary  \ref{intersection}, and  assume in addition that $M$ is relative Cohen--Macaulay with respect to $P \cap Q$. Then $\rdim(Q, M)=\cd(P, M)=1$. The converse holds when $M$ is relative Cohen--Macaulay with respect to $P$, as well.
\end{Corollary}
\begin{proof}
Let $M$ be a relative Cohen--Macaulay with respect to $P \cap Q$. Corollary \ref{intersection} yields
$ \cd(P \cap Q, M)=\dim M-1=\grade(P \cap Q, M)$.
Note that $H^{i}_{I\cap J}(-)=H^{i}_{IJ}(-)$ for all $i$ and for all graded ideals $I$ and $J$ of any graded Noetherian ring $R$.
Hence
\[
\grade(P \cap Q, M)=\grade(PQ, M)=\min\{\grade (P, M), \grade (Q, M)\}.
\]
The second equality follows from \cite [Proposition 9.1.3(b)]{BH}. Thus we conclude that
$\dim M-1=\min\{\grade (P, M), \grade (Q, M)\}$. Here we consider two cases: Let  $\grade (P, M)\leq \grade (Q, M)$.
 Then Theorem \ref{rdim1} yields
 \[
 \dim M-1=\grade (P, M)\leq \grade (Q, M)=\rdim (Q, M)=\dim M-\cd(P, M).
 \]
 Thus $\cd(P, M)=1$.  We also see that $\rdim(Q, M)=1$, because
 \[
 0<\rdim (Q, M)=\dim M-1=\grade (P, M)\leq \cd(P, M)=1.
 \]
  Now we assume that $\grade (Q, M)\leq \grade (P, M)$. Similarly we obtain $\rdim(Q, M)=\cd(P, M)=1$. The converse follows from Theorem \ref{rdim1}.
 \end{proof}
 \section{Maximal relative Cohen--Macaulay}
 \begin{Definition}{\em
 Let $M$ be a finitely generated bigraded $S$-module. We call $M$ to be maximal relative Cohen--Macaulay with respect to $Q$ if $M$ is relative Cohen--Macaulay with respect to $Q$ such that $\rdim (Q, M)=\dim K[y]=n$. We also call $M$ to be maximal relative Cohen--Macaulay with respect to $P$ if $M$ is relative Cohen--Macaulay with respect to $P$ such that $\rdim (P, M)=\dim K[x]=m$.}
 \end{Definition}
 \begin{Remark}{\em
 In Lemma \ref{example} we observe that if $J=0$, then the ring $R$ is maximal relative Cohen--Macaulay with respect to $Q$ with  $\rdim (Q, R)=\dim K[y]=n$ and if $I=0$,  then $R$ is maximal relative Cohen--Macaulay with respect to $P$ with  $\rdim (P, R)=\dim K[x]=m$. In Example \ref{2} the ring $R$
 is relative Cohen--Macaulay with respect to $Q$ but not maximal relative Cohen--Macaulay with respect to $Q$. }
 \end{Remark}
 \begin{Proposition}
  Let $M$ be a finitely generated bigraded $S$-module. Then the following hold:
 \begin{itemize}
\item[{(a)}]  $M$ is maximal relative Cohen--Macaulay with respect to $Q$ if and only if $y_1, \dots, y_n$  is an $M$-sequence.
\item[{(b)}]  $M$ is maximal relative Cohen--Macaulay with respect to $P$ and $Q$ if and only if $M$ is a free $S$-module.
 \end{itemize}
 \end{Proposition}
 \begin{proof}
 Part (a) follows from \cite[Corollary 1.6.19]{BH} and the fact that  $\cd (Q, M)\leq n$. For the prove (b) since $M$ is relative Cohen--Macaulay with respect to $P$ and $Q$, then Theorem \ref{rdim1}  implies that $\dim M=\dim S$ and  Corollary \ref{Cohen--Macaulay} yields that $M$ is Cohen--Macaulay. Therefore $M$ is maximal Cohen--Macaulay $S$-module and hence by the Auslander-Buchsbaum formula it follows that $M$ is a free $S$-module. The other implication is obvious.
 \end{proof}

\begin{center}
{Acknowledgment}
\end{center}
\hspace*{\parindent}
I would like to thank Professor J\"urgen Herzog for his helpful comments.

\end{document}